\documentclass[11pt]{amsart}

\newtheorem{theorem}{Theorem}[section]
\newtheorem{lemma}[theorem]{Lemma}

\theoremstyle{definition}
\newtheorem{definition}[theorem]{Definition}

\newtheorem{example}[theorem]{Example}

\newtheorem{proposition}[theorem]{Proposition}
\newtheorem{corollary}[theorem]{Corollary}
\newtheorem{remark}[theorem]{Remark}

\theoremstyle{remark}

\newcommand{\be}{\begin{equation}}
\newcommand{\ee}{\end{equation}}

\numberwithin{equation}{section}

\begin{document}

\title{Circle action, lower bound of fixed points and characteristic numbers}

\author{Ping Li}
\address{Department of Mathematics, Tongji University, Shanghai, 200092, China}
\address{Department of Mathematics and Information Sciences, Tokyo
Metropolitan University, Tokyo 192-0397, Japan}

\email{pingli@tongji.edu.cn}
\thanks{The author was partially
supported by National Natural Science Foundation of China (Grant No.
11101308) and JSPS Postdoctoral Fellowship for Foreign Researchers.}




\dedicatory{Dedicated to Professor Richard Palais on the occasion of
his 80th birthday}


\begin{abstract}
Given an $S^1$-manifold with isolated fixed points, some recent
papers are concerned with the relationship between the least number
of fixed points and the characteristic numbers of this manifold, and
their proofs have some similar features. The main purpose of this
short survey article is, by using the language of equivariant
cohomology, to present a unified method to deal with such problems,
of which the related known results are direct corollaries.
\end{abstract}

\maketitle

\section{Introduction}
In transformation group theory, it is a classical and important
topic to study various topological obstructions (such as the
vanishing or non-vanishing of certain characteristic numbers) to the
existence of non-trivial circle actions on manifolds with specified
properties.

Recently, some papers (\cite{FR}, \cite{LL1}, \cite{LL2}, \cite{PT}) are concerned with the least number of fixed points of
circle action on manifolds, and their proofs share some similarity. The main purpose of this note is to unify
the results in the above-mentioned papers into the framework of equivariant cohomology, which allows us to treat various results simultaneously.

Our paper is organized as follows. In Section \ref{section2}, after
introducing some necessary symbols and notation related to
equivariant cohomology, we will prove some properties based on the
Localization Theorem. As applications of these properties,we give
some examples in Section \ref{section3}, which are main results of
(\cite{FR}, \cite{LL1}, \cite{LL2}, \cite{PT}) and can be derived
from our properties in Section \ref{section2}. In addition to these
examples, some related remarks are also given in this section.

Throughout the following two sections, all $S^1$-manifolds mentioned
are connected, smooth, closed and oriented and all circle actions on
such manifolds are smooth. When we say an $S^1$-almost-complex
manifold $(M^{2n},J)$, we mean that the manifold is closed and given
the canonical orientation and the circle action preserves the
almost-complex structure $J$. We use superscripts to denote the
\emph{real} dimensions of the manifolds.

\section{Main observations}\label{section2}
As usual, we use $S^1$ to denote the circle group. All cohomology groups will have $\mathbb{Q}$-coefficients unless otherwise stated. Given an $S^1$-space $X$, its $S^1$-equivariant cohomology is written as $H^{\ast}_{S^1}(X):=H^{\ast}(\textrm{E}S^1\times_{S^1}X)$, where $\textrm{E}S^1\times_{S^1}X$ is the standard Borel construction. $H^{\ast}_{S^1}(X)$ can be viewed as a module over the coefficient ring $H^{\ast}_{S^1}(\{\textrm{pt}\})=H^{\ast}(\textrm{B}S^1)=\mathbb{Q}[z]$, where $z\in H^2(\textrm{B}S^1)$ is the Euler class
of the universal complex line bundle over $\textrm{B}S^1$. The inclusion $i:~\{1\}\rightarrow S^1$ induces a restriction map $i^{\ast}:~H^{\ast}_{S^1}(X)\rightarrow H^{\ast}(X)$.

Let $M^{2m}$ be an $S^1$-manifold and $M^{S^1}$ be the fixed point set of this action. For an equivariant cohomology class $v\in H^{2k}_{S^1}(M)$ and a fixed point $x\in M^{S^1}$, we write the restriction of $v$ to $H^{2k}_{S^1}(\{x\})=\mathbb{Q}z^k$ as $v_xz^k$, $v_x\in\mathbb{Q}$. The tangent bundle $\textrm{T}M$ of $M$ can be lifted to an oriented vector bundle $\textrm{E}S^1\times_{S^1}\textrm{T}M$ over $\textrm{E}S^1\times_{S^1}M$ and we denote its Euler class as $e\in H^{2m}_{S^1}(M)$.

If $M^{S^1}$ is finite, we have the following Localization Formula (\cite{AB}), which is largely due to Bott.

\begin{theorem}[Localization Formula]\label{localization}
Suppose that $M^{S^1}$ is finite and $w\in H^{2m}_{S^1}(M)$. Then we have
\be\label{localizationexpression}i^{\ast}(w)[M]=\sum_{x\in M^{S^1}}\frac{w_x}{e_x}\in\mathbb{Q}.\ee
\end{theorem}

\begin{remark}\label{fixedfreeremark}
If $M^{S^1}$ is empty, i.e., $M$ admits a fixed point free $S^1$-action, this Localization Formula tells us that $i^{\ast}(w)[M]=0$ for all $w\in H^{2m}_{S^1}(M)$.
\end{remark}

Note that when $x\in M^{S^1}$ is an isolated fixed point, $e_x$ is nonzero and thus the right-hand side of (\ref{localizationexpression}) makes sense.

\begin{proposition}\label{prop}
Let $u\in H^{2k}_{S^1}(M)$ and $v\in H^{2l}_{S^1}(M)$ such that $m=kr+l$. If $i^{\ast}(u^rv)[M]\neq 0$, then $\sharp M^{S^1}\geq r+1$. Here $\sharp M^{S^1}\geq r+1$ is the cardinality of $M^{S^1}$.
\end{proposition}
\begin{proof}
Applying (\ref{localizationexpression}) to the elements $w=u^s\cdot v\cdot z^{k(r-s)}$ ($0\leq s\leq r$), we have

\begin{eqnarray}\label{loca2}
\sum_{x\in M^{s^1}}\frac{u_x^s\cdot v_x}{e_x}
=\left\{\begin{array}{l}0,~\textrm{if}~0\leq s<r\\
i^{\ast}(u^rv)[M],~\textrm{if}~s=r.
\end{array}
\right.\end{eqnarray}
Since $i^{\ast}(u^rv)[M]\neq 0$, we know that the following set

\be\label{set}\{u_x~|~x\in M^{S^1},~v_x\neq 0\}\subset\mathbb{Q}\ee

is nonempty. So we can list the elements in (\ref{set}) in increasing order as
$a_1<a_2<\cdots<a_n$ $(n\geq 1)$. Let $\lambda_i$ $(1\leq i\leq n)$ be the sum of the terms $\frac{v_x}{e_x}$ such that $v_x\neq 0$ and $u_x=a_i$. So from (\ref{loca2}) we have

\begin{eqnarray}\label{loca3}
\sum_{i=1}^na_i^s\cdot\lambda_i
=\left\{\begin{array}{l}0,~\textrm{if}~0\leq s<r\\
i^{\ast}(u^rv)[M],~\textrm{if}~s=r.
\end{array}
\right.\end{eqnarray}

We assert that $n\geq r+1$. In fact, if $n\leq r$, the $n$-order matrix $(a_i^j)_{1\leq i,j\leq n}$
is nonsingular and we deduce that all $\lambda_i=0$ and so $i^{\ast}(u^rv)[M]=0$, a contradiction. Therefore $n\geq r+1$. But by definition $\sharp M^{S^1}\geq n$, which gives the desired proof.
\end{proof}

\begin{remark}
The trick used in this Proposition can be traced back to Pelayo-Tolman
(\cite{PT}, Lemma 8) and the present author and Liu (\cite{LL1}, Lemmas
3.1 and 3.2). L\"{u}-Tan also
used it in their proof of Theorem 1.1 in \cite{LT}.
\end{remark}

From the proof of Proposition \ref{prop} we have the following corollary.

\begin{corollary}\label{lemma}
Let the notation be as above. We have
\begin{enumerate}
\item
if $n\leq r$, then all $\lambda_i=0$. In particular, $i^{\ast}(u^rv)[M]=0$.

\item
if $n\leq r+1$ and $i^{\ast}(u^rv)[M]=0$, then all $\lambda_i=0$.
\end{enumerate}
\end{corollary}

For the critical case $\sharp M^{S^1}=r+1$, we have the following corollary.

\begin{corollary}\label{coro}
Suppose that $\sharp M^{S^1}=r+1$. Then
\begin{enumerate}
\item
if $n<\sharp M^{S^1}=r+1$, then all $\lambda_i=0$. In particular, $i^{\ast}(u^rv)[M]=0$.

\item
if $n=\sharp M^{S^1}=r+1$, then $i^{\ast}(u^rv)[M]\neq 0$.
\end{enumerate}
\end{corollary}

\begin{proof}
(1) is directly from Item (1) of Corollary \ref{lemma}. For (2),
if $i^{\ast}(u^rv)[M]=0$, from Item (2) of Corollary \ref{lemma}
 we deduce that all $\lambda_i=0$, which means for every $i$,
  the number of fixed points $x$ in $M^{S^1}$ such that
$u_x=a_i$ and $\frac{v_x}{e_x}=\lambda_i$ is at least two. Thus
$n\leq\frac{1}{2}\sharp M^{S^1}$, which contradicts to the
assumption.
\end{proof}

\section{Applications}\label{section3}
In this section, we first recall some basic knowledge on equivariant
vector bundles and particularly on $S^1$-complex line bundles. Then
using these, together with the properties proved in Section
\ref{section2}, we give some examples, which are main results in
some related previous papers.

\subsection{Preliminaries}
Given any real (resp. complex) $S^1$-vector bunde $\xi$ (resp. $\eta$) over $M$, we have the corresponding equivariant Pontrjagin classes (resp. Chern classes) $p_j^{S^1}(\xi)=p_j(\textrm{E}S^1\times_{S^1}\xi)\in H^{4j}_{S^1}(M)$ \big(resp. $c_j^{S^1}(\eta)=c_j(\textrm{E}S^1\times_{S^1}\eta)\in H^{2j}_{S^1}(M)$\big) which lift $p_j(\xi)\in H^{4j}(M)$ \big(resp. $c_j(\eta)\in H^{2j}(M)\big)$ in the sense that $i^{\ast}(p_j^{S^1}(\xi))=p_j(\xi)$ \big(resp. $i^{\ast}(c_j^{S^1}(\eta))=c_j(\eta)\big).$ In particular, the tangent bundle $\textrm{T}M$ is a real $S^1$-vector bundle over $M$ such that $i^{\ast}(p_j^{S^1}(M))=p_j(M)$. If $(M^{2m},J)$ is an almost-complex manifold, the holomorphic tangent bundle $\textrm{T}^{\prime}M$, which is in the sense of $J$, is a complex $S^1$-vector bundle such that $i^{\ast}(c_j^{S^1}(M))=c_j(M)$.

Suppose $M$ is an $S^1$-manifold and $c\in H^{2}(M;\mathbb{Z})$. It
is well known that there exists a complex line bundle $L$ over $M$,
which is unique up to bundle isomorphism, such that $c_{1}(L)=c$. We
call $c$ \emph{admissible} if the given $S^1$-action on $M$ can be
lifted to an action on the corresponding complex line bundle $L$. It
is well-known that (\cite{HY}, \cite{Mu}) $c$ is admissible if and
only if $c\in i^{\ast}(H^{2}_{S^1}(M,\mathbb{Z}))$. The following
very useful lemma is quite well-known (cf. \cite{HY}, \cite{Mu}).
But for the convenience of the readers, we provide a short proof.

\begin{lemma}\label{lemma2}
If the first Betti number $b_1(M)=0$, then any element in $ H^{2}(M;\mathbb{Z})$ is admissible.
\end{lemma}
\begin{proof}
First we note that if $b_1(M)=0$, then $H^1(M,\mathbb{Z})=0$ by the
Universal Coefficient Theorem.

The map $i^{\ast}:~H^{\ast}_{S^1}(M,\mathbb{Z})\rightarrow
H^{\ast}(M,\mathbb{Z})$ fits into a long exact sequence

\be 0\rightarrow
H^1_{S^1}(M,\mathbb{Z})\xrightarrow{i^{\ast}}H^1(M,\mathbb{Z})\rightarrow
H^0_{S^1}(M,\mathbb{Z})\xrightarrow{\cdot
z}H^2_{S^1}(M,\mathbb{Z})\xrightarrow{i^{\ast}}H^2(M,\mathbb{Z})\rightarrow
H^1_{S^1}(M,\mathbb{Z})\rightarrow\cdots,\nonumber\ee from which we
deduce that $H^1_{S^1}(M,\mathbb{Z})=0$ and thus
$i^{\ast}:~H^{2}_{S^1}(M,\mathbb{Z})\rightarrow H^{2}(M,\mathbb{Z})$
is surjective.
\end{proof}

\subsection{Examples}
The following example is a generalization of a result of Fang-Rong (\cite{FR}, Theorem 1.1).

\begin{example}
Let $M^{2m}$ be a manifold with $b_1(M)=0$. If $M$
admits a fixed point free $S^1$-action, then for all $c\in
H^{2}(M;\mathbb{Z})$ and $0\leq j\leq [\frac{m}{2}]$, we have
$$(c^{m-2j}\cdot L_j)[M]=0,$$
where $L_j$ is any polynomial of degree $j$ in the subalgebra of $H^{\ast}(M)$ generated by Pontrjagin classes $p_i(M)$ ($\textrm{deg}(p_i)=i$). In particular, $(c^{n})[M]=0$.
\end{example}
\begin{proof}
By Lemma \ref{lemma2} there exists $\bar{c}\in H^2_{S^1}(M)$ such that $i^{\ast}(\bar{c})=c$. Then in Proposition \ref{prop}, we can take $k=1$, $r=m-2j$ and $v$ to be any polynomial of degree $j$ in the equivariant Pontrjagin classes $p^{S^1}_{\ast}(M)$.
\end{proof}

\begin{remark}
The orginal statement of (\cite{FR}, Theorem 1.1) is only for those $L_j$ in the
Hirzebruch $L$-polynomial of $M$: $L=L_{0}+L_{1}+\cdots+\L_{[\frac{m}{2}]}$.
\end{remark}

The following example is given by the author and Liu (\cite{LL1}, \cite{LL2}).

\begin{example}\label{ll}
Suppose $(M^{2m},J)$ (resp. $N^{4m}$) is an $S^1$-almost-complex manifold (resp. smooth $S^1$-manifold). Let $(1^{\lambda_1}2^{\lambda_2}\cdots m^{\lambda_m})$ be a partition of $m$, which means $\lambda_1+2\lambda_2+\cdots +m\lambda_m=m$ and all $\lambda_i$ are nonnegative integers. If the corresponding Chern number
$c_1^{\lambda_1}\cdots c_m^{\lambda_m}[M]\neq 0$ (resp. Pontrjagin number $p_1^{\lambda_1}\cdots p_m^{\lambda_m}[N]\neq 0$), then
the circle action on $M$ (resp. $N$) has at least $\textrm{max}\{\lambda_1,\ldots,\lambda_m\}+1$ fixed points.
\end{example}
\begin{proof}
Note that, under the above condition, $\sharp M^{S^1}$ (resp. $\sharp N^{S^1}$) is nonempty. Thus we may assume that $\sharp M^{S^1}$ (resp. $\sharp N^{S^1}$) is finite. Hence we can take $k=j$ (resp. $k=2j$) ($1\leq j\leq m$), $r=\lambda_j$, $u=c_j^{S^1}(M)$ (resp. $u=p_j^{S^1}(N)$) and
$v=\prod_{i\neq j}[c_i^{S^1}(M)]^{\lambda_i}$ (resp. $v=\prod_{i\neq j}[p_i^{S^1}(N)]^{\lambda_i}$) in Proposition \ref{prop} and let $j$ range over $1,\ldots,m$.
\end{proof}

\begin{remark}Define

$$\lambda_1:=\textrm{max}\{\textrm{max}\{\lambda_1,\ldots,\lambda_m\}~|~c_1^{\lambda_1}\cdots c_m^{\lambda_m}[M]\neq 0\}+1$$
and
$$\lambda_2:=\textrm{max}\{\textrm{max}\{\lambda_1,\ldots,\lambda_m\}~|~p_1^{\lambda_1}\cdots p_m^{\lambda_m}[N]\neq 0\}+1$$
respectively.
Then Example \ref{ll} in fact shows that the circle action on $M$ (resp. $N$) has at least $\lambda_1$ (resp. $\lambda_2$) fixed points.
\end{remark}

Before giving another two examples, we need to introduce some more notation.

Let $M^{2m}$ be an $S^1$-manifold and $c\in H^{2}(M;\mathbb{Z})$ be
admissible such that $c_1(L)=c$. We fix a lifting of the action on
$L$. Moreover, if the given $S^1$-action has isolated fixed points,
say $P_1,\cdots,P_{\mu}$, we consider the fiber
$L_{p_{i}}\cong\mathbb{C}$ of $L$ over the fixed point $P_i$. Each
fiber $L_{p_{i}}$ is an $S^1$-module, and hence corresponds to an
integer $a_i$ such that the action of each element $\lambda\in S^1$
on $L_{p_{i}}$ is given by multiplying $\lambda^{a_i}$. The integers
$a_i$ will be called the weights of $L$ (with respect to the given
lifting) at the fixed point $P_i$. Note that if we choose another
lifting of this $S^1$-action then the weights $a_i$ are changed
simultaneously to $a_{i}+a$ for some integer $a$. So we give the
following definition.

\begin{definition}Suppose $c$ is admissible with respect to an $S^1$-action having isolated
fixed points. We call $c$ \emph{everywhere injective} if the weights
$a_1,\cdots,a_{\mu}$ are mutually distinct. Similarly, we call such $c$
\emph{somewhere injective} if these $a_1,\cdots,a_{\mu}$ have the
property that, at least one of these integers is different from the
other $\mu-1$ integers.
\end{definition}

Clearly the property of everywhere injective implies that of
somewhere injective.

The following example is a generalization of some much earlier results of Hattori (\cite{Ha1},
p.441-p.443), which was drawn my attention by Kefeng Liu.

\begin{example}\label{main result}
Suppose $M^{2m}$ is a smooth, oriented $S^{1}$-manifold and $c\in H^{2}(M;\mathbb{Z})$ is
admissible.
\begin{enumerate}
\item
If $(c)^m[M]\neq 0$, then the given action has at least $m+1$ fixed
points.

\item
If the given action has isolated fixed points and $c$ is somewhere
injective, then the number of the fixed points is at least
$m+1$.

\item
Suppose the given circle action has exactly $m+1$ fixed points. Then
there are exactly two possibilities: \begin{enumerate}
\item
$c$ is everywhere injective and $(c)^m[M]\neq 0$;
\item $c$ is not somewhere injective and $(c)^m[M]=0$.
\end{enumerate}
\end{enumerate}
\end{example}

\begin{proof}
Suppose $i^{\ast}(\bar{c})=c$ where $\bar{c}\in H^2_{S^1}(M)$. Then we can take
in Proposition \ref{prop} $k=1$, $r=m$, $u=\bar{c}$ and $v=1$. In this case
$u_{P_i}=a_i$ and $v_{P_i}=1$.

Now (1) follows directly from Proposition \ref{prop}. If $c$ is
somewhere injective, then some $\lambda_i$ must be nonzero. Hence
(2) follows from Item (1) of Corollary \ref{lemma}. (3) follows from
Corollary \ref{coro}.
\end{proof}

\begin{remark}
The original purpose of \cite{Ha1} was to investigate the
constraints on certain circle actions which have similar features to
those arising from the positivity of curvature, as mentioned by its
author in the Introduction. But \cite{Ha1} contains many interesting
results which are of independent interest and deserves further
attention.
\end{remark}

A special case of Item (2) in Example \ref{main result} is due to Pelayo-Tolman (\cite{PT}, Theorem 1).
\begin{example}\label{PT}
Suppose $(M^{2m},J)$ is an $S^1$-almost-complex manifold with isolated fixed points. If $c_1(M)$ is somewhere injective, then the number of the fixed points is at least
$m+1$.
\end{example}

\section*{Acknowledgments} I would like to thank Professor Kefeng Liu for drawing my attention to the reference
\cite{Ha1}, and Professor Paul Baum for many fruitful discussions related to this paper during his visit to Fudan
University in Shanghai. Special thanks should go to Professor Michael Crabb, who read the original version of this paper carefully, shared his many new insights with me and helped me to formulate this paper as the current form.

\end{document}